\documentclass[11pt]{article}
\usepackage[margin=1in]{geometry}
\usepackage[T1]{fontenc}
\usepackage{lmodern,microtype}
\DisableLigatures[-]{encoding = T1, family = tt*}
\usepackage{amsmath,amssymb,amsthm,mathtools}
\usepackage{booktabs,longtable,array}
\usepackage{todonotes}

\usepackage{needspace}
\usepackage{xcolor,tikz}
\usetikzlibrary{arrows.meta,positioning,calc}
\usepackage{hyperref}
\hypersetup{colorlinks=true,linkcolor=blue,citecolor=blue,urlcolor=blue,
pdftitle={Hyperbolic origami surfaces of arbitrary genus},
pdfsubject={Explicit eight-vertex blocks and a computer-assisted construction}}
\newtheorem{theorem}{Theorem}[section]
\newtheorem{lemma}[theorem]{Lemma}
\newtheorem{proposition}[theorem]{Proposition}

\theoremstyle{definition}\newtheorem{definition}[theorem]{Definition}
\theoremstyle{remark}
\newcommand{\HH}{\mathbb H}
\newcommand{\RR}{\mathbb R}

\newcommand{\TT}{\mathcal T}
\newcommand{\CC}{\mathcal C}
\newcommand{\MM}{\mathcal M}
\newcommand{\LL}{\mathcal L}

\newcommand{\Int}{\operatorname{int}}
\newcommand{\sgn}{\operatorname{sgn}}

\newcommand{\code}[1]{\texttt{\detokenize{#1}}}

\newcommand{\ClaimRootRadius}{ 10^{-35} }

\title{Hyperbolic origami surfaces of arbitrary genus}
\author{Aly Rajwani}
\date{\today}
\begin{document}
\maketitle
\begin{abstract}
  For any $g \geq 2$, we construct a hyperbolic origami surface of genus $g$.
\end{abstract}

\tableofcontents

\section{Introduction}\label{sec:introduction}
\subsection{Main Result}\label{subsec:main-result}

Recently, Schwartz constructed an eight-vertex Euclidean paper torus \cite[Theorem~1.2]{Schwartz}, and Zou introduced its hyperbolic analogue and constructed a ten-vertex example of genus two \cite[Theorem~1.1]{Zou}. Using Schwartz's triangulation and Zou's robust-embedding approach, we construct hyperbolic origami surfaces of every genus $g \geq 2$.

For a closed hyperbolic surface $\Sigma$ and a finite geodesic triangulation $\TT$ of $\Sigma$, a \emph{hyperbolic origami surface} is an embedding $S:(\Sigma, \TT) \hookrightarrow \HH^3$ such that $S$ restricts to an isometry on each triangle. Concretely, it is a collection of hyperbolic triangles glued in $\HH^3$ such that the cone angle at each vertex is $2\pi$. Here is our main theorem: 

\begin{theorem}\label{thm:main}
  For every $g \geq 2$, there exists a closed, orientable, genus $g$ hyperbolic surface $S_g: (\Sigma_g, \TT_g) \hookrightarrow \HH^3$ with $5g + 3$ vertices, $21g + 3$ edges, and $14g + 2$ faces. 
\end{theorem}

\begin{figure}[ht]
    \centering
    \begin{minipage}{0.48\textwidth}
        \centering
        \includegraphics[width=\textwidth]{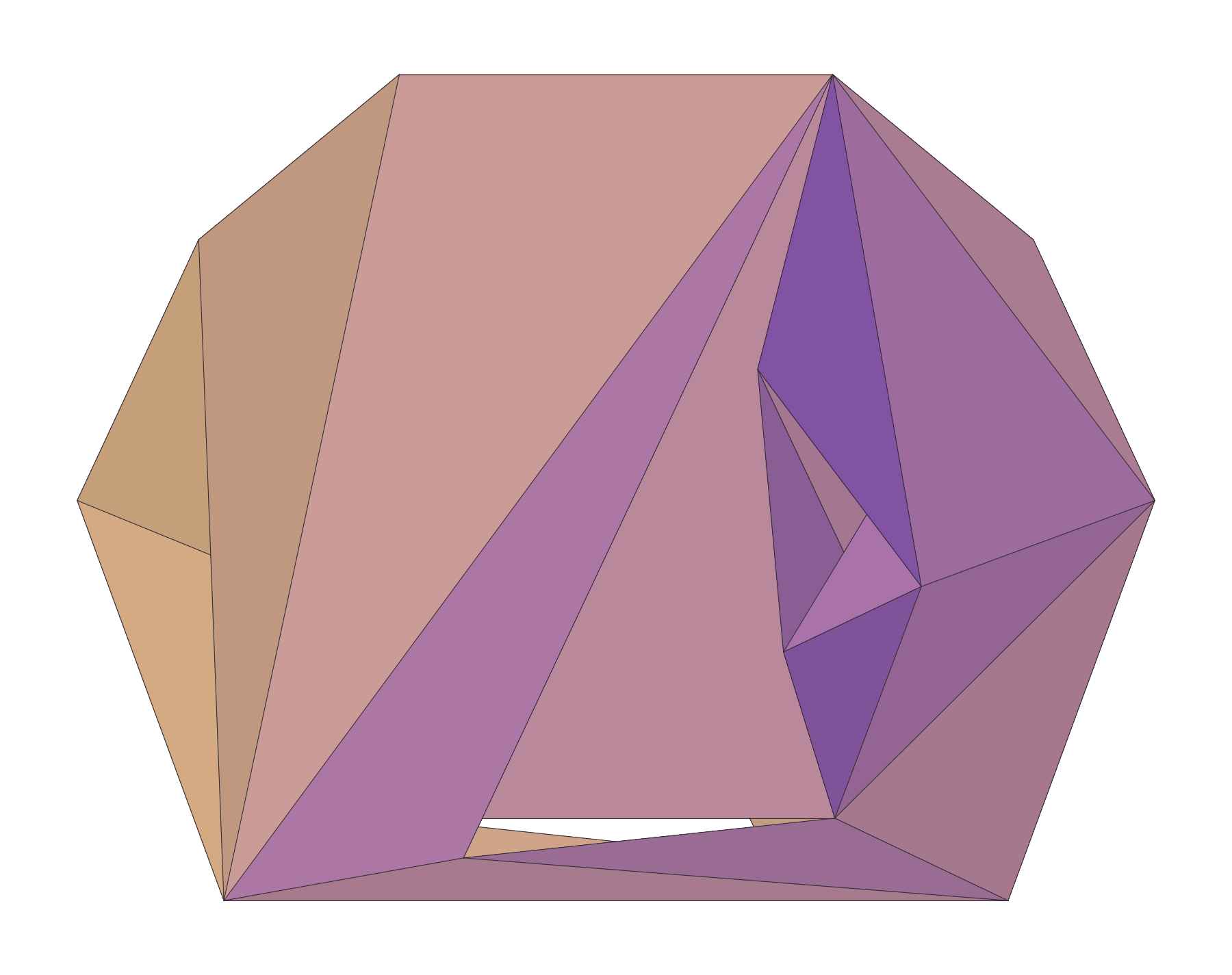}
    \end{minipage}
    \hfill
    \begin{minipage}{0.48\textwidth}
        \centering
        \includegraphics[width=\textwidth]{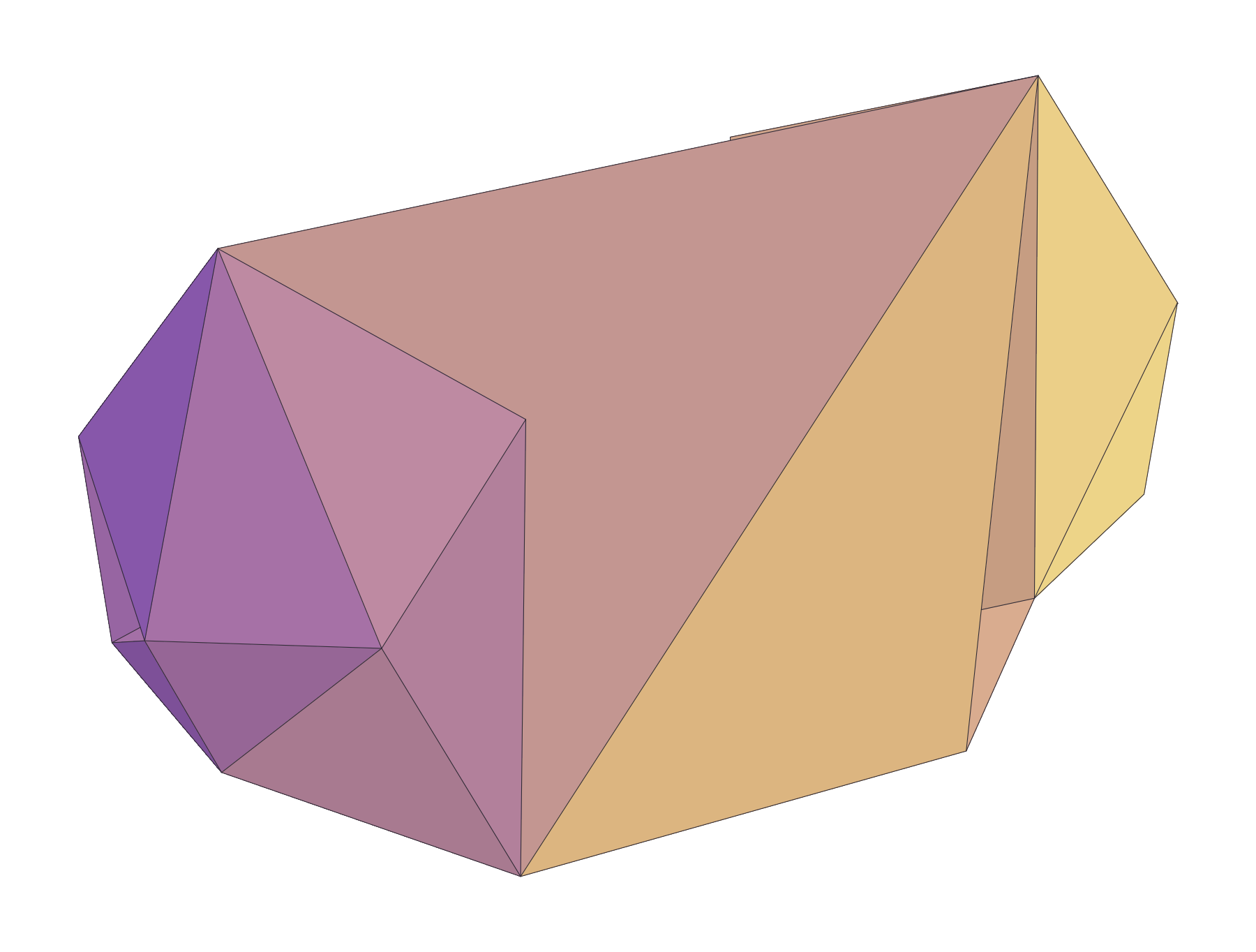}
    \end{minipage}
    \caption{Two views a genus $3$ hyperbolic origami surface in the Beltrami--Klein model.}
    \label{fig:comparison}
\end{figure}

We construct $S_g$ by joining two caps to a chain of middle pieces, each with eight vertices. We first prove that the cap and middle piece have the required cone angles at each vertex, then prove that we can assemble them without self-intersection. Note that this does not produce vertex-minimal triangulations; Zou's ten-vertex example is smaller than the 13-vertex example here. Figure \ref{fig:comparison} depicts a genus-three example.

\subsection{Accompanying Programs}

Readers can access the accompanying code at the following link:
\begin{center}
  \url{https://github.com/alyrajwani/origami-surfaces}
\end{center}

\subsection{AI Statement}

GPT-6.1 Astra was used in the development and implementation of the code.

\subsection{Acknowledgements}

I am incredibly grateful to Rich Schwartz and Zhengyu (Byron) Zou who inspired me to work on this problem. The methodology of this paper closely mirrors \cite{Schwartz} and \cite{Zou}, and this paper could not have been written without their work, advice, and input.

\section{The Blocks and Their Geometry}\label{sec:geometry}

\subsection{The Beltrami--Klein Model}\label{subsec:bk-model}

Let $\HH^3$ be the simply connected Riemannian three-manifold of constant sectional curvature $-1$. In this paper, we work primarily in the Beltrami--Klein model of $\HH^3$, which has point-set 
 \[\HH^3=B_1(0)=\{X=(x_1,x_2,x_3)\in\RR^3:\|X\|<1\}.\]
For $X \in \HH^3$, write $a_X = 1 - \|X\|^2$. The Riemannian metric and its inner product are
\begin{equation}\label{eq:kleinmetric}
  g|_X = \frac{\sum_{i=1}^3(dx_i)^2}{a_X} + \frac{(\sum_{i=1}^3x_i\,dx_i)^2}{a_X^2},\qquad 
  \langle V,W\rangle_X =\frac{a_X\langle V,W\rangle+\langle X,V\rangle\langle X,W\rangle}{a_X^2}.
\end{equation}
Here, $\langle \cdot, \cdot\rangle$ and $\|\cdot\|$ denote the Euclidean inner product and norm respectively. Write $\|V\|_X=\sqrt{\langle V,V\rangle_X}$. Our reason for working in this model is that geodesics are straight lines, and totally geodesic planes are affine planes intersected with the ball $B_1(0)$. A triangular face is the Euclidean convex hull of its three vertices as a point set, with the metric \eqref{eq:kleinmetric}. Thus, we often perform computations using Euclidean coordinates. We sometimes write $\mathbb B^3=B_1(0)$.

Vertices are denoted $X_i$, corner angles at $X_i$ in the face $(X_i, X_j, X_K)$ by $\theta_{ijk}$, and their sum, i.e., the total corner angle at $X_i$, by $\theta_i$. We hope that corner angles for interior vertices are $2\pi$ and for boundary vertices in the cap and middle pieces to be $\pi$, which we encode by the vector 
\[\Theta(\TT) = (\theta_0 - \tau_0, \dots, \theta_7 - \tau_7),\]
where $\tau_i$ is $2\pi$ for interior vertices and $\pi$ for boundary vertices, and the indices range from $0$ to $7$ since both pieces are based on Schwartz's eight-vertex triangulation.  

\subsection{The Candidate Models}\label{subsec:blocks}

We build the surface from two fixed pieces. The \emph{cap} $C$ is a torus with one disk removed, and the \emph{middle} $M$ is a torus with two disks removed. Their boundary curves are triangles, which we call \emph{ports}, and we call the planes they determine \emph{port planes}. At each interior vertex the target angle is $2\pi$, and at each boundary vertex the target is $\pi$ so that two pieces glued along a port contribute $2\pi$ at each seam vertex.

The cap lies on one side of its port plane. The middle lies between two ultraparallel port planes, i.e., the planes have positive distance. Thus, reflections across these planes produce a chain of middles with two caps closing the ends such that the resulting surface is embedded. Each added middle adds one handle, so two caps and $g-2$ middles give genus $g$.

All labels in the following paragraph are taken modulo eight. Let $\TT$ have vertices $0,\ldots,7$ and oriented faces
\begin{equation}\label{eq:faces}
  f_i=(i,i+1,i+3),\qquad h_i=(i,i+3,i+2),\qquad 0\leq i\leq7.
\end{equation}
This triangulation is nice in the sense that its degree sequence is $(6, 6, 6, 6, 6, 6, 6, 6)$. This is the triangulation used by Schwartz in \cite{Schwartz}.

\subsubsection{The Cap}\label{subsubsec:cap}

The cap complex is $\CC = \TT \setminus \Int(f_2)$. It retains all vertices and edges, and has $15$ faces. Its boundary is the cycle $(3,2,5)$, written in the orientation induced from the remaining faces. For the spatial realization $C$, all three boundary vertices have $z=0$, while every interior vertex has $z>0$. The desired retained-angle conditions are
\begin{equation}\label{eq:captarget}
 \theta_i(C) = \begin{cases}\pi,&i\in\{2,3,5\},\\
   2\pi,&i\in\{0, 1, 4, 6, 7\}.
 \end{cases}
\end{equation}
Here $\theta_i$ sums only the angles of the retained faces. 

\subsubsection{The Middle}\label{subsubsec:middle}

The middle complex is $\MM = \TT \setminus (\Int(f_5)\cup\Int(h_7))$. It retains all vertices and edges and has 14 faces. Its boundary consists of cycles $(6, 5, 0)$ and $(2, 7, 1)$. Note that the removed triangles share no vertices and that the permutation $j(i) = 7 - i$ preserves the face set and exchanges the ports. This is the vertex pairing of Schwartz's pup tent. We impose the ambient symmetry
\begin{equation}\label{eq:J}
  J(x,y,z) = (-x,-y,z),\qquad X_{7-i} = JX_i, \qquad 0 \leq i \leq 3.
\end{equation}
It is a hyperbolic isometry because it is an orthogonal transformation of the Klein ball. The desired angle sums are $\pi$ at the six boundary vertices and $2\pi$ at vertices $3$ and $4$. By imposing the symmetry $J$, it suffices to impose the four equations at $0,1,2,3$.

\subsection{The Rational Candidates}\label{subsec:rational-candidates}

Write $\widehat C$ and $\widehat M$ for the two candidates with rational coordinates. Their exact coordinates are listed in Appendix~\ref{app:centers}, and Figure~\ref{fig:blocks} shows their realizations. The reader can verify the coordinate tables against those files by running \path{code/proofs/check_manuscript.py}. Section~\ref{sec:root} constructs exact blocks near these candidates by solving the angle and matching equations.

\begin{figure}[ht]
    \centering
    \begin{minipage}{0.40\textwidth}
        \centering
        \includegraphics[width=\textwidth]{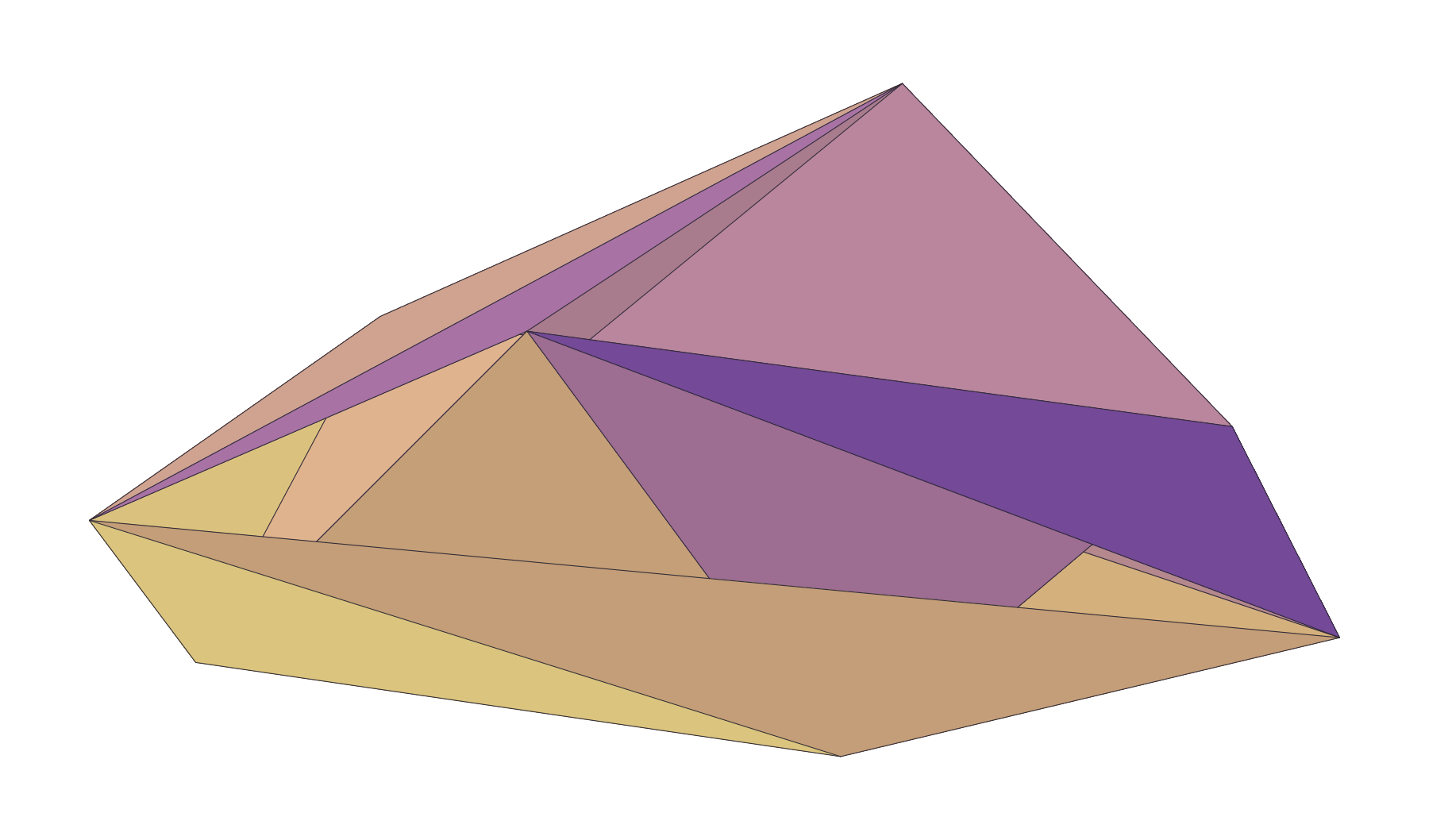}
    \end{minipage}
    \hfill
    \begin{minipage}{0.40\textwidth}
        \centering
        \includegraphics[width=\textwidth]{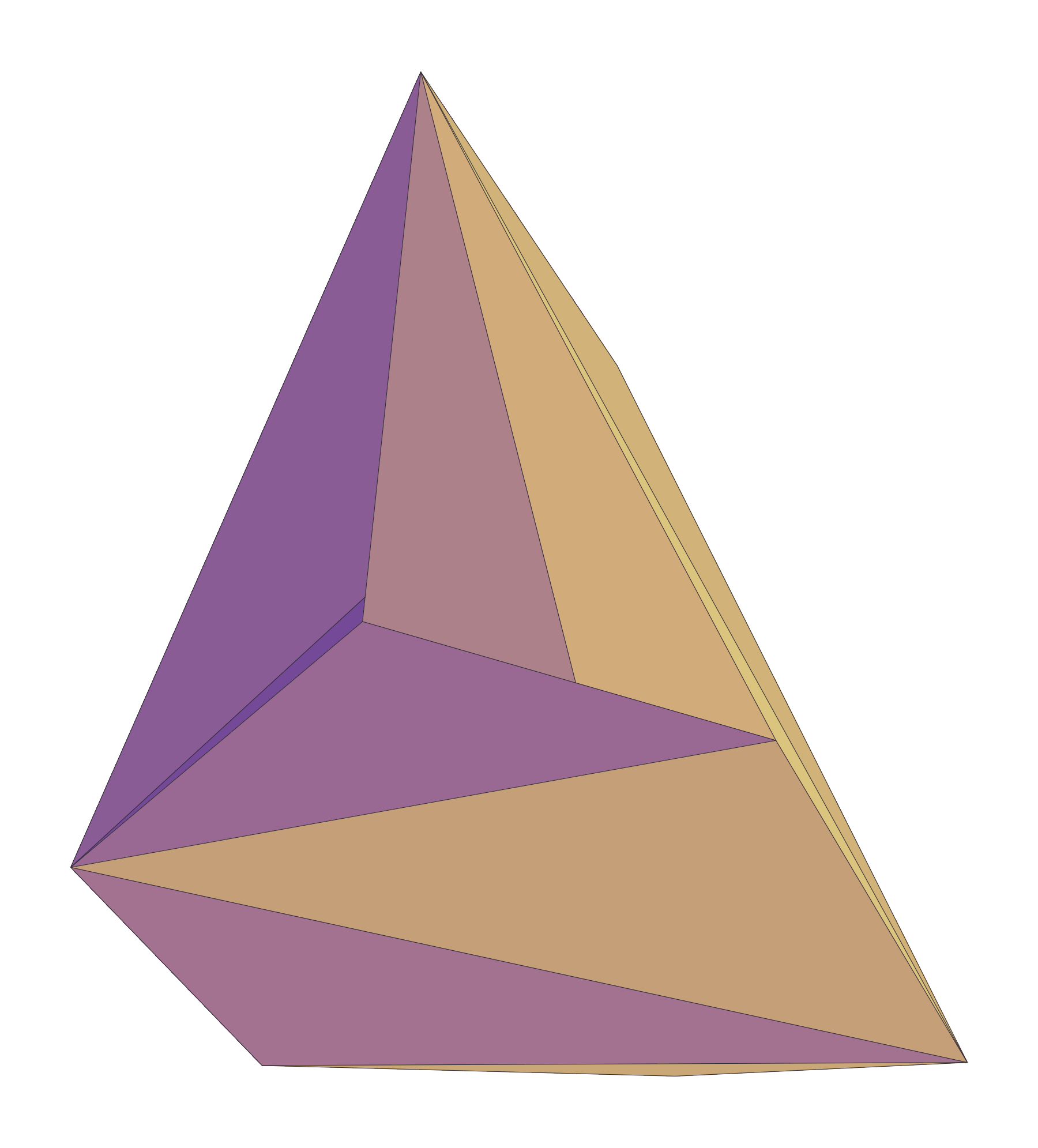}
    \end{minipage}
    \caption{The rational candidates. On the left is the cap and on the right is the middle.}
    \label{fig:blocks}
\end{figure}

\subsection{Corner Angles}\label{subsec:corner-angles}

For a triangle $(X, Y, Z)$, set $V = Y - X$ and $W = Z - X$, and define 
\[B_{VW} = a_X\langle V, W\rangle + \langle X, V\rangle\langle X, W\rangle, \qquad B_V = B_{VV},\]
and consider the functions 
\begin{equation}\label{eq:A}
  A(X, Y, Z) = \frac{\langle V, W\rangle_X^2}{\langle V, V\rangle_X \cdot \langle W, W\rangle_X} = \frac{B_{VW}^2}{B_V \cdot B_W}, \qquad \sigma = \sgn\langle V, W\rangle_X = \sgn B_{VW}.
\end{equation}
Then the corner angle at $X$ in the triangle is 
\begin{equation}\label{eq:corner-angle}
  \theta = \arccos\left(\sigma\sqrt{A(X,Y,Z)}\right) = \operatorname{atan2}\left(\sqrt{B_VB_W - B_{VW}^2}, B_{VW}\right).
\end{equation}
Note that for a nondegenerate triangle, $B_V, B_W > 0$ and $B_VB_W - B_{VW}^2 > 0$. 

\subsection{A Common Port Shape}

For distinct $p,q\in\mathbb B^3$ define
\begin{equation}\label{eq:inv}
 I(p,q) = \frac{(1-p\cdot q)^2}{(1-\|p\|^2)(1-\|q\|^2)} = \cosh^2d_{\HH^3}(p,q).
\end{equation}
For rational Klein coordinates this invariant is rational. On positive lengths, $\ell\mapsto\cosh^2\ell$ is strictly increasing. Thus, we check equality of ports by checking equality after applying $I$. This provides an exact boundary matching condition without introducing transcendental length values. Specifically, three rational constants $\kappa_0,\kappa_1,\kappa_2$ in Appendix~\ref{app:invariants} are fixed once and for all. They were computed using \eqref{eq:inv} from an initial rational middle candidate. We impose the conditions
\begin{align}
  I(C_3,C_2)&=I(M_5,M_6)=\kappa_0,\notag\\
  I(C_2,C_5)&=I(M_6,M_0)=\kappa_1,\label{eq:matching}\\
  I(C_5,C_3)&=I(M_0,M_5)=\kappa_2.\notag
\end{align}
Since the two middle ports are congruent by \eqref{eq:J}, the cap matches either one.

\subsection{Interval Arithmetic}\label{subsec:interval-arithmetic}

Much of this paper depends on precise computer arithmetic. In order to perform this, we enclose quantities by an interval containing it with rational endpoints and center and perform calculations using these rational coordinates. Using the Python package \code{mpmath}, we can do arithmetic to arbitrary precision, which, for this paper, was 100 decimals. Note that since the interval arithmetic $\code{mpmath}$ in experimental for some functions, calculations for inverse trigonometric functions were done using the file \code{rigorous_angles.py}. See Appendix \ref{app:angles} provides detail on the validity of our calculations.

\section{The Candidate Blocks}\label{sec:candidate-blocks}

In this section, we establish that our candidates have cone angles and port planes quite close to those specified in sections \ref{subsubsec:cap} and \ref{subsubsec:middle}, and that they remain embedded under perturbation. The strategy closely follows Schwartz and Zou. 

\subsection{Robust Flatness}

A \emph{hyperbolic polyhedral surface} is an embedding $S: (\Sigma, \TT) \hookrightarrow \HH^3$ such that the image of a triangle in $\TT$ under $S$ is a geodesic triangle. Thus, a hyperbolic origami surface is a hyperbolic polyhedral surface whose cone angles are all $2\pi$. 

\begin{definition}\label{def:flat}
A polyhedral realization of either block is \emph{$\epsilon$-flat} if
\[\max_i|\theta_i-\tau_i|<\epsilon,\qquad
  \tau_i = \begin{cases}
    \pi,& i\text{ is a boundary vertex},\\
    2\pi,& i\text{ is an interior vertex}.
  \end{cases}\]
\end{definition}

\begin{proposition}\label{prop:flat}
The candidates $\widehat C$ and $\widehat M$ are $10^{-80}$-flat.
Further, each of their three ordered port invariants differs from its prescribed
value $\kappa_j$ by less than $10^{-80}$.
\end{proposition}

\begin{proof}
  The proof is via computation. We evaluate equations \eqref{eq:corner-angle} and \eqref{eq:inv} at each vertex of both candidates using the methods described in section \ref{subsec:interval-arithmetic}, and find that the difference from the desired values all lie in $(-3 \cdot 10^{-89}, 3 \cdot 10^{-89})$. The skeptical reader may verify the computations by running the following command:
  \begin{center}
    \begin{verbatim}
    python3 -E -B code/proofs/gluing_certify_matched_blocks.py --verify-only
    \end{verbatim}
  \end{center}
\end{proof}

\subsection{Robust Embeddedness}

We next verify that these candidates remain embedded under small changes in their coordinates. Here, we take advantage of the Beltrami--Klein model by forgetting the hyperbolic structure and checking embeddedness using Euclidean geometry. 

\begin{definition}\label{def:robust}
  An embedded polyhedral block $S$ with vertices $X_0,\ldots,X_7$ is \emph{$\lambda$-robustly embedded} if every realization with the same faces and vertices $X'_i$ satisfying
  \[\max_i\|X'_i-X_i\|_\infty\leq\lambda\]
  is an embedding with nondegenerate triangles.
\end{definition}

\begin{lemma}\label{lem:pairs}
  Let $S$ be a realization of either block by nondegenerate triangles. If each pair of distinct closed faces intersects exactly in its common abstract simplex, then $S$ is embedded.
\end{lemma}

\begin{proof}
  See \cite[Lemma 3.2]{Schwartz} or \cite[Lemma 3.5]{Zou}.
\end{proof}

\begin{lemma}\label{lem:disjoint}
  Let $T_1$ and $T_2$ be faces with disjoint vertex sets. Suppose an affine functional $\ell(x)=n\cdot x-b$ satisfies
  \[\min_{p\in T_1^{(0)}}\ell(p)\geq\mu>0, \qquad \max_{q\in T_2^{(0)}}\ell(q)\leq-\mu,\]
where $T_i^{(0)}$ denotes the vertices of $T_i$. If each vertex changes by at most $\delta$ in infinity norm and $\delta\|n\|_1<\mu$, the perturbed faces remain disjoint.
\end{lemma}

\begin{proof}
  If $T_i'$ is obtained from $T_i$ by perturbing the vertices by at most $\delta$, the value of each functional changes by at most $\delta \|n\|_1$, which is less than $\mu$ by assumption. Thus, $\ell$ witnesses a hyperplane separating the vertices of $T_1'$ from $T_2'$, and by convexity and linearity of the dot product, this implies $T_1'$ and $T_2'$ are disjoint. 
\end{proof}

\begin{lemma}\label{lem:vertex}
  Let $T_1=(v,p_1,p_2)$ and $T_2=(v,q_1,q_2)$. Suppose an affine functional $\ell(x)=n\cdot x-b$ satisfies
  \[\min_{j=1,2}\bigl(\ell(p_j)-\ell(v)\bigr)\geq\mu>0, \qquad \max_{j=1,2}\bigl(\ell(q_j)-\ell(v)\bigr)\leq-\mu.\]
  If every vertex changes by at most $\delta$ in coordinate infinity norm and $2\delta\|n\|_1<\mu$, the perturbed faces intersect only at their common vertex.
\end{lemma}

\begin{proof}
  Let $T_i'$ be obtained from $T_i$ by perturbing $(v, p_1, p_2, q_1, q_2)$ by at most $\delta$. Since $\ell(p)-\ell(v)=n\cdot(p-v)$, each difference changes by at most $2\delta\|n\|_1 < \mu$, so 
  \[\min_{j = 1, 2} (\ell(p'_j) - \ell(v')) > 0, \qquad \max_{j = 1, 2} (\ell(q'_j) - \ell(v')) < 0.\]
  Suppose $T_1'$ and $T_2'$ intersect at a point $x$, which we may express as 
  \[x = t_0v' + t_1p_1' + t_2p_2' = s_0v' + s_1q_1' + s_2q_2'\] 
  where $t_0 + t_1 + t_2 = s_0 + s_1 + s_2 = 1$ and $t_i, s_j \geq 0.$ It follows that 
  \[t_1(\ell(p'_1) - \ell(v')) + t_2(\ell(p_2') - \ell(v')) = s_1(\ell(q_1') - \ell(v')) + s_2(\ell(q_2') - \ell(v')),\]
  but our previous inequality implies this is only possible if $t_1 = t_2 = s_1 = s_2 = 0$, which yields $x = v'$.
\end{proof}

\begin{lemma}\label{lem:edge}
  Let $T_1=(a,b,c)$ and $T_2=(a,b,d)$ be nondegenerate triangles, and set
  \[\ell_{abc}(x)=\bigl((b-a)\times(c-a)\bigr)\cdot(x-a).\]
  If $\ell_{abc}(d)=\det(b-a,c-a,d-a)\ne0$, then the triangles intersect precisely in their common edge.
\end{lemma}

\begin{proof}
  The nonvanishing is $\ell_{abc}(d) = \det(b - a, c - a, d - a)$ is precisely the statement that $T_1$ and $T_2$ lie in different planes, hence they can only intersect in their common edge. 
\end{proof}

\begin{proposition}\label{prop:robust}
  The candidates $\widehat C$ and $\widehat M$ are $10^{-8}$-robustly embedded.
\end{proposition}
\begin{proof}
  Let $\widehat X_0, \dots, \widehat X_7$ denote the vertices of $\widehat C$ or $\widehat M$. For each of the 24 coordinate entries we use the rational interval $[\widehat X_i^k - 10^{-8}, \widehat X_i^k + 10^{-8}]$. Interval arithmetic yields the following strict lower bounds throughout the boxes:
  \begin{center}\small
    \begin{tabular}{lcc}\toprule
      Quantity & Cap & Middle\\\midrule
      $1-\|X_i\|^2$ & $1/3$ & $1/5$\\
      $\|(X_b-X_a)\times(X_c-X_a)\|^2$ & $1/20$ & $1/8$\\
      Signed shared-edge determinant & $1/1000$ & $1/100$\\
      Signed separating functional value & $99/100$ & $99/100$\\\bottomrule
    \end{tabular}
  \end{center}
  The first row demonstrates that our perturbations all lie in $\HH^3$ and the second that all triangles are nondegenerate. The third and fourth rows provide values of $\mu$ which allow us to apply Lemmas \ref{lem:disjoint}--\ref{lem:edge} and conclude that distinct faces in the described perturbations intersect in their abstract simplex. Lemma \ref{lem:pairs} them implies all such perturbations are embedded, hence $\widehat C$ and $\widehat M$ are $10^{-8}$-robustly embedded. The skeptical reader can verify these numerical bounds by running \path{code/proofs/verify_all_coordinate_robustness.py}.
\end{proof}

\subsection{Port Planes}\label{sec:support}

In proving result about our port planes, we restrict to the family of embeddings specified in Section \ref{sec:systems}. Each parameter varies by at most $r = \ClaimRootRadius$ from its rational center, and we impose the cap plane being $\{z = 0\}$ and the order two symmetry of the middle.  

\begin{proposition}\label{prop:placement}
  Throughout this family, the cap meets $z=0$ exactly in its port boundary and lies in $z>0$ otherwise. The middle meets each of two ultraparallel supporting planes exactly in the corresponding port boundary and lies otherwise strictly between the planes. Every port triangle is nondegenerate.
\end{proposition}
\begin{proof}
  The cap boundary vertices satisfy $z = 0$ by construction, and interval arithmetic yields $z > 0$ for the other five vertices throughout the family. By convexity and our use of the Beltrami--Klein model, the cap meets $z = 0$ exactly in its port boundary. Order the ports for the middle piece as $(5, 6, 0)$ and $(7, 2, 1)$ so that the outward normal for each port $(a, b, c)$ is 
  \[n = (X_b - X_a) \times (X_c - X_a),\qquad d = n \cdot X_a.\]
  Interval arithmetic yields $\|n\|^2 > 0$ and $d - n \cdot X_i > 0$ for every vertex not in the port, including vertices in the other port. Convexity implies that the middle lies between the port planes and meets these planes only in the corresponding port boundary. We now prove these planes are ultraparallel using the Lorentz model of hyperbolic space. Let $\langle x, y\rangle_\mathcal{L} = -x_0y_0 + x_1y_1 + x_2y_2 + x_3y_3$ and let $\mathcal{H} = \{x \in \RR^{1, 3} : \langle x, x\rangle_\mathcal{L} = -1, x_0 > 0\}$. Write $N_i = (d_i, n_i) \in \RR^{1, 3}$. Interval arithmetic yields 
  \begin{equation}\label{eq:ultra}
    q_i := \langle N_i, N_i\rangle_\mathcal{L} > 0, \qquad L := \langle N_0, N_1\rangle_\mathcal{L} < 0, \qquad L^2 > q_0q_1. 
  \end{equation}
  Let $V = \operatorname{span}(N_0, N_1)$, so the Gram matrix of the Lorentz inner product restricted to $V$ is 
  \[G = \begin{pmatrix} q_0 & L \\ L & q_1 \end{pmatrix}\]
    has negative determinant, hence has eigenvalues with opposite signs. Since the signature is $(-, +, +, +)$, this implies that $V^\perp$ is positive definite. Suppose the planes intersect at a point $x$, so $\langle x, N_0\rangle_\mathcal{L} = \langle x, N_1\rangle_\mathcal{L} = 0$. This implies that $x \in V^\perp$, but since $x \in \mathcal{H}$, we also have $\langle x, x\rangle_\mathcal{L} = -1$, a contradiction. This also proves that the planes cannot meet at an ideal point, i.e., a point satisfying $\langle x, x\rangle_\mathcal{L} = 0$, so the planes are ultraparallel. The skeptical reader can verify these computations by reading Appendix \ref{sec:dyadic} and running the files \path{code/proofs/gluing_certify_matched_blocks.py} and \path{code/exact_geometry/check_certificates.py}.
\end{proof}

\section{Finding Origami Surfaces}\label{sec:root}

Thus far, our work has been approximate. Following the methods of Schwartz and Zou, we demonstrate that we can improve these results to exact origami surfaces.

\subsection{The Equations}\label{sec:systems}

Write $C_i = (x_i, y_i, z_i)$ for cap coordinates and $M_i$ for middle coordinates. We fix all coordinates except for the following 11 which are free: 
\begin{equation}\label{eq:capvars}
  (x_5, x_3, y_6, y_1, x_2, y_4, y_3, y_0, z_7, z_0, x_1).
\end{equation}
For the middle, the free coordinates are 
\begin{equation}\label{eq:midvars}
  (z_0, x_1, y_2, x_3, z_1, x_2, x_0).
\end{equation}
The coordinates for $M_4, \dots, M_7$ are determined by the imposed symmetry. See Appendix \ref{app:centers} for exact values of the fixed coordinates for both the cap and middle. Before selecting free coordinates, the Jacobians at the rational centers have sizes $11 \times 21$ for the cap and $7 \times 12$ for the middle. QR decomposition with column-pivoting selects the coordinates listed in \eqref{eq:capvars} and \eqref{eq:midvars}, giving square matrices of sizes $11 \times 11$ and $7 \times 7$, respectively. For $B = C, M$, let $c_B$ be the vector of free coordinates at the rational centers and define $Q_B = \{u : \|u - c_B\|_\infty \leq r\}$. Let $U_B$ be the open set of free parameters for which all vertices lie in the open Klein ball and all faces are nondegenerate. The bounds in \ref{sec:dyadic} give the inclusion $Q_B \subset U_B$. We compute corner angles using \eqref{eq:corner-angle}. We define a cap map $F_C: U_C \to \RR^{11}$ coordinate-wise: its first eight coordinates are 
\[\Theta_i = \theta_i - \tau_i, \quad i = 0, \dots, 7, \quad \tau_i = \begin{cases} \pi,& i \text{ is a boundary vertex,} \\ 2\pi, & i \text{ is an interior vertex}, \end{cases}\]
  and its last three are $I(p_j, p_{j+1}) - \kappa_j$ with indices modulo 3, where $(p_0, p_1, p_2)$ is a port and the expression comes from the port invariants specified in \eqref{eq:matching}. The middle map $F_M: U_M \to \RR^7$ has first four coordinates 
\[(\Theta_0, \Theta_1, \Theta_2, \Theta_3) = (\theta_0 - \pi, \theta_1 - \pi, \theta_2 - \pi, \theta_3 - 2\pi),\]
and its last three are also $I(p_j, p_{j+1}) - \kappa_j$ as specified in \eqref{eq:matching}. All denominators and radicands remain strictly positive throughout the boxes, so both $F_C$ and $F_M$ are smooth, which we use in Lemma \ref{lem:contraction}. For each block, we have a rational square matrix $R_B$ supplied in the code, and we define $H_B = R_BF_B$. These matrices were obtained by approximating the inverse of the Jacobian $DF_B$.  

\subsection{Finding a Zero}

We find an exact zero to our approximate solution by constructing a function to which we can apply Banach's fixed point theorem.

\begin{lemma}\label{lem:contraction}
    Let $r > 0$ and let $F: U \to \RR^m$ be continuously differentiable on an open neighborhood $U \subset \RR^m$ of the closed box $B = \{x \in \RR^m : \|x - c\|_\infty \leq r\}$. Let $R$ be a square $m \times m$ matrix, and suppose 
  \[\|RF(c)\|_\infty \leq \eta, \quad \sup_{x \in B} \|I - RDF(x)\|_\infty \leq q < 1, \quad \frac{\eta}{1 - q} < r.\]
  Then $F$ has a unique zero $z$ in $B$, and it satisfies $\|z - c\|_\infty \leq \eta/(1-q)$. 
\end{lemma}
\begin{proof}
  Let $G(x) = x - RF(x)$, so $DG = I - RDF$, and by our assumptions, $\|DG(x)\|_\infty \leq q$ in $B$. By mean value estimates, we have $\|G(x) - G(y)\|_\infty \leq q\|x - y\|_\infty$ for all $x, y \in B$. Furthermore, since $G(c) - c = - RF(c)$, we have for all $x \in B$ that 
  \begin{align*}
    \|G(x) - c\|_\infty &\leq \|G(x) - G(c)\|_\infty + \|G(c) - c\|_\infty \\
    & \leq q\|x - c\|_\infty + \eta \\
    &\leq qr + \eta \\
    &< r 
  \end{align*}
  hence $G$ is a contraction mapping of $B$. By Banach's fixed point theorem, $G$ has a unique fixed point $z$. Standard linear algebra shows that $\|I - RDF(x)\|_\infty < 1$ implies $RDF(x)$ is invertible, and as these are square matrices, $R$ is invertible. Thus, we have the implications 
  \[G(z) = z \Longleftrightarrow RF(z) = 0 \Longleftrightarrow F(z) = 0,\]
  so $F$ has a unique zero in $B$, and the bound $\|z - c\|_\infty \leq \eta/(1 - q)$ follows from our previous computation. 
\end{proof}

\subsection{The Jacobian}

Write $u = (u_1, \dots, u_m)$, where $m = 11$ for the cap and $m = 7$ for the middle, in the orders \eqref{eq:capvars} and \eqref{eq:midvars}. Recall that for a face $(i, j, k)$, the angle at $i$ is denoted $\theta_{ijk}$, so for an angle coordinate, we have 
\[\frac{\partial \Theta_i}{\partial u_\ell} = \sum_{(i, j, k) \text{ a face}} \frac{\partial \theta_{ijk}}{\partial u_\ell}.\]
Consider a triangle with vertices $X, Y, Z$. Write $V = Y - X, W = Z - X, a = 1 - \|X\|^2, \alpha = \langle X, V\rangle,$ and $\beta = \langle X, W\rangle$, where the norm and inner product here are Euclidean. In this notation, the quantities from subsection \ref{subsec:corner-angles} are 
\[B_V = a\langle V, V\rangle + \alpha^2, \quad B_W = a\langle W, W\rangle + \beta^2, \quad B_{VW} = a\langle V, W\rangle + \alpha\beta.\]
Let $d$ denote some linear combination of $\partial_X, \partial_Y, \partial_Z$. The product rule yields 
\begin{align*}
  dB_V &= \langle V, V\rangle da + 2a\langle V, dV\rangle + 2\alpha d\alpha, \\
  dB_W &= \langle W, W\rangle da + 2a\langle W, dW\rangle + 2\beta d\beta, \\
  dB_{VW} &= \langle V, W\rangle da + a\langle dV, W\rangle + a\langle V, dW\rangle + \alpha d\beta + \beta d\alpha.
\end{align*}
Now write $s = \sqrt{B_VB_W - B_{VW}^2}$ and $\theta = \operatorname{atan2}(s, B_{VW})$, so 
\[\cos \theta = \frac{B_{VW}}{\sqrt{B_VB_W}}, \qquad \sin \theta = \frac{s}{\sqrt{B_VB_W}}.\]
Differentiating the first expression yields 
\[-\sin\theta d\theta = \frac{dB_{VW} - (dB_V/B_V + dB_W/B_W) \cdot B_{VW}/2}{\sqrt{B_VB_W}},\]
and substituting the second yields 
\begin{equation}\label{eq:dangle}
  d\theta = \frac{B_{VW}/2 \cdot (dB_V/B_V + dB_W/B_W) - dB_{VW}}{\sqrt{B_VB_W - B^2_{VW}}}.
\end{equation}
Summing the corner-angle derivatives at each vertex gives the first eight rows of $DF_C$ and the first four rows of $DF_M$. The last three rows of each Jacobian come from differentiating the port invariant in \eqref{eq:inv}. For endpoints $p, q$, write $b = 1 - \langle p, q\rangle, a_p = 1 - \langle p, p\rangle,$ and $a_q = 1 - \langle q, q\rangle$ so that $I(p, q) = b^2/(a_pa_q)$. Differentiating yields 
\[dI = I(p, q)\left(-\frac{2(\langle p, dq\rangle + \langle dp, q\rangle)}{b} + \frac{2\langle p, dp\rangle}{a_p} + \frac{2\langle q, dq\rangle}{a_q}\right),\]
and since the target $\kappa_j$ are fixed, $d(I - \kappa_j) = dI$. The verifier evaluates these derivatives by defining pairs $(v, \dot{v})$, where $v$ is an interval and $\dot{v}$ is a vector of intervals enclosing its gradient, and using automatic differentiation. Evaluation over the coordinate box $Q_B$ encloses every Jacobian entry at every point, so multiplication by the exact rational matrix $R_B$ gives an enclosure of $I - R_BDF_B$ throughout the box. We obtained the following bounds on all the quantities in Lemma \ref{lem:contraction}, and the skeptical reader can verify the following bounds with the file \code{code/proofs/gluing_certify_matched_blocks.py}:

\begin{center}
  \begin{tabular}{lll}\toprule
    Quantity & Cap & Middle\\\midrule
    Box radius $r$ & $10^{-35}$ & $10^{-35}$\\
    $\|F(c)\|_\infty$ & $<10^{-80}$ & $<10^{-80}$\\
    $\eta=\|RF(c)\|_\infty$ & $<10^{-60}$ & $<10^{-60}$\\
    $q=\sup_{Q_B}\|I-RDF\|_\infty$ & $<10^{-20}$ & $<10^{-20}$\\
  \bottomrule\end{tabular}
\end{center}

\begin{proposition}\label{prop:blocks}
  There exist embedded realizations $C, M$ of $\CC, \MM$ by nondegenerate geodesic triangles with corner angle $2\pi$ at interior vertices, corner angle $\pi$ at boundary vertices, congruent ports, and with the port plane properties of Proposition \ref{prop:placement}.  
\end{proposition}
\begin{proof}
  The bounds in the above table yield 
  \[\eta + qr < 10^{-60} + 10^{-20}\ClaimRootRadius < \ClaimRootRadius = r,\]
  so the hypotheses of Lemma \ref{lem:contraction} are satisfied, and $F_B$ has a zero in $Q_B$. The angle coordinates of $F_C, F_M$ give the required angle conditions, and the port conditions give equality of the three corresponding side lengths, hence congruence of the ordered port triangles. Propositions \ref{prop:robust} and \ref{prop:placement} apply throughout $Q_B$, whence the result follows. 
\end{proof}

\section{Proof of the Main Theorem}\label{sec:chain}

\begin{proof}[Proof of Theorem \ref{thm:main}] 
  Fix $g \geq 2$ and let $n = g - 2$. Start with the cap $C$ from Proposition \ref{prop:blocks}. Its port lies in a plane $P$, and the rest of the cap lies entirely in one side of $P$. If $n > 0$, we glue the first port of $M$ to the port of the cap, which is possible because of their congruence. 

  If $n > 0$, match the first port of $M$ to that of the cap. Proposition \ref{prop:blocks} gives congruence of the ports, so matching their vertices matches each geodesic edge. At corresponding vertices, take an orthonormal tangent frame matched by the plane isometry. Completing them by unit normals gives an isometry carrying one frame to the other, so we can place $M$ across $P$ from the cap. At each later step, reflect the most recently glued middle across its free port plane to obtain a chain with one additional piece separated by a plane. Let $Q$ be the other port plane of the new middle. By Proposition \ref{prop:placement}, $Q$ is disjoint from $P$ and lies entirely on the new side of $P$, so we may repeat this process arbitrarily. 

  After gluing $n$ middle pieces, glue the second cap to the free port and place it on the exterior side by the earlier choice of normal. The only intersections are the prescribed boundary cycles, so the induced map from the compact glued surface to $\HH^3$ is a continuous injection, hence an embedding. When $n = 0$, this construction simply reflects the first cap across its port plane. 

  The surface is closed, connected, and orientable by construction. Since we glue two caps and $n$ middles along circles, the Euler characteristic is $\chi = 2(-1) + n(-2)$, and since it is closed, the genus is $(2-\chi)/2 = n + 2 = g$. Each vertex not on a seam has has a cone angle of $2\pi$, and each vertex on a seam has a cone angle of $\pi + \pi = 2\pi$. The incident sectors therefore form a hyperbolic disk, which gives the glued surface its hyperbolic metric with each face realized isometrically.

  Each of he $n + 1$ seams identifies three pairs of vertices and three pairs of edges, so we have a triangulation with 
  \[(V, E, F) = 2(8, 24, 15) + n(8, 24, 14) - (n+1)(3, 3, 0) = (5g + 3, 21g + 3, 14g + 2).\]
\end{proof}

\begin{figure}[ht]
    \centering
    \includegraphics[width=0.4\textwidth]{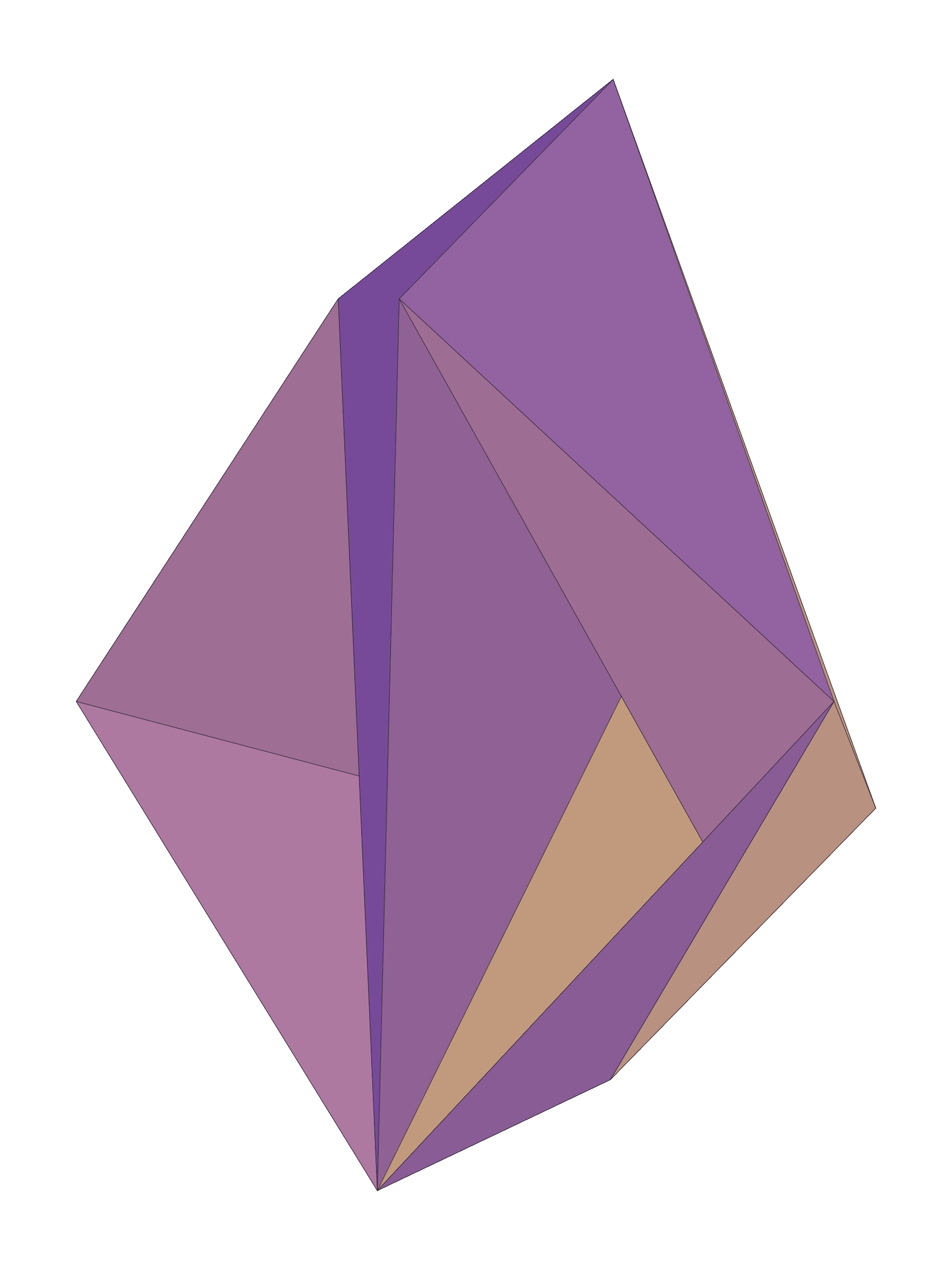}
    \caption{Two caps give a genus 2 surface.}
    \label{fig:first}
\end{figure}

\begin{figure}[ht]
    \centering
    \includegraphics[width=0.5\textwidth]{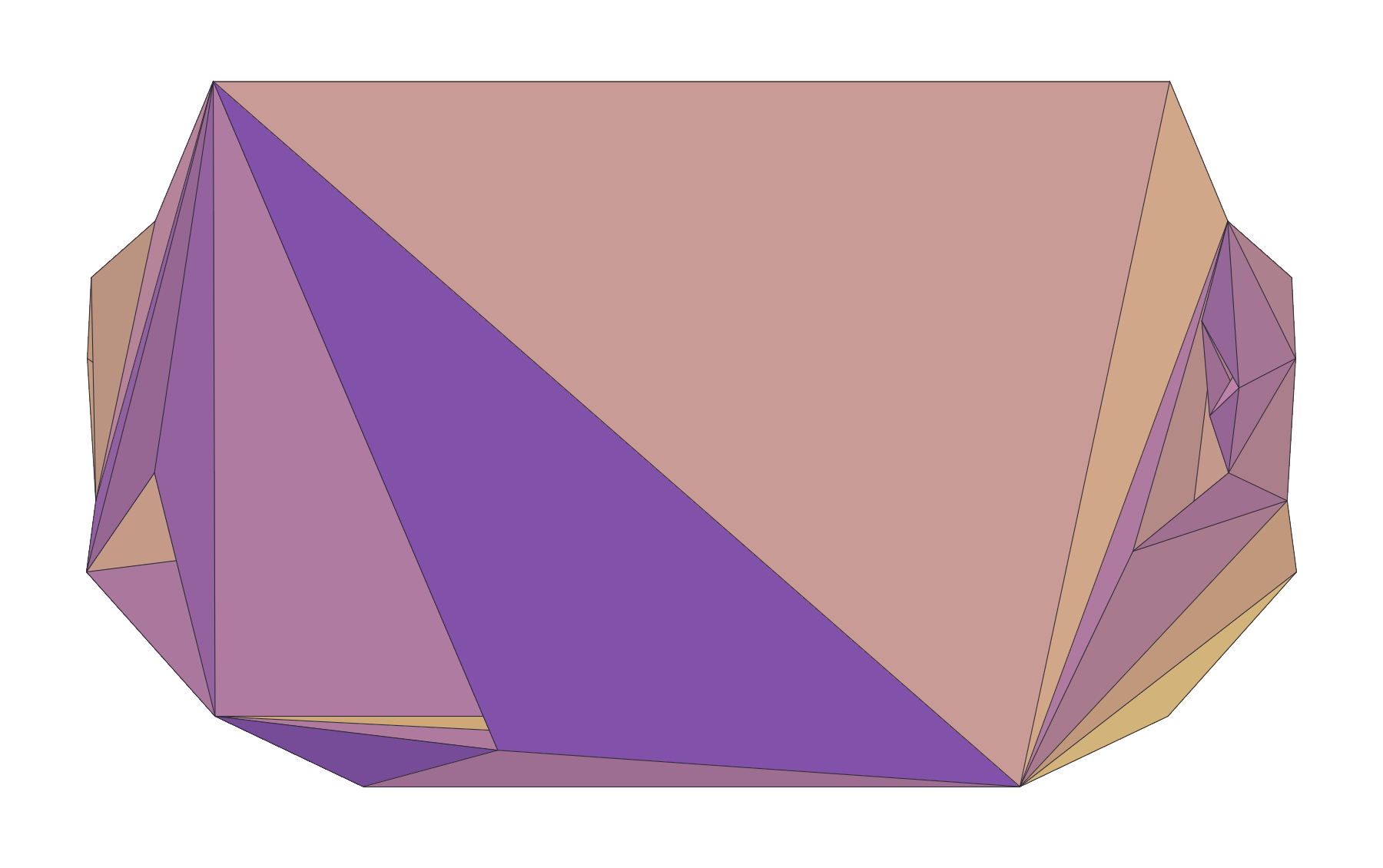}
    \caption{Two caps and three middles give a genus 5 surface.}
    \label{fig:second}
\end{figure}
\newpage
\appendix
\section{Exact Centers and Fixed Coordinate Parameters}\label{app:centers}

In the following table, a finite decimal denotes an exact rational number. An asterisk in a fixed-parameter table indicates a variable coordinate. The long decimals for variable centers are split across lines, and their concatenation denotes the exact decimal value. 

\subsection{Cap Fixed Parameters}
\begin{center}\scriptsize
  \begin{tabular}{rrrr}\toprule
    $i$ & $x_i$ & $y_i$ & $z_i$\\\midrule
    0 & -0.026858804929243224 & $\ast$ & $\ast$ \\
    1 & $\ast$ & $\ast$ & 0.03670350228395235 \\
    2 & $\ast$ & 0.15756525598542975 & 0 \\
    3 & $\ast$ & $\ast$ & 0 \\
    4 & 0.4221185460436074 & $\ast$ & 0.014557430710943203 \\
    5 & $\ast$ & -0.1956870501260478 & 0 \\
    6 & 0.2831137865258996 & $\ast$ & 0.025684907531724616 \\
    7 & -0.17639985786323786 & 0.38554490406720143 & $\ast$ \\
  \bottomrule\end{tabular}
\end{center}

\Needspace{12\baselineskip}
\subsection{Cap Variable Centers}
\small
\begin{longtable}{rrcl}
  \toprule Index & Vertex & Coordinate & Exact rational center \\\midrule
  \endfirsthead
  \toprule Index & Vertex & Coordinate & Exact rational center \\\midrule
  \endhead
  1 & 5 & x & \shortstack[l]{\texttt{0.6800372983574048664080031598}\\\texttt{920330629179356881737151561811}\\\texttt{857134673490596889726965696742}\\\texttt{34}} \\[5pt]
  2 & 3 & x & \shortstack[l]{\texttt{-0.559180811381750871126054864}\\\texttt{229769707397853292278480904030}\\\texttt{729783230114895592230331162189}\\\texttt{947}} \\[5pt]
  3 & 6 & y & \shortstack[l]{\texttt{-0.651895447935056208716360788}\\\texttt{996991227509512042572370560065}\\\texttt{078013110744756797295565813380}\\\texttt{262}} \\[5pt]
  4 & 1 & y & \shortstack[l]{\texttt{0.5747158845771109455038408213}\\\texttt{108799805996383440951491744055}\\\texttt{493977396913379231722365669558}\\\texttt{51}} \\[5pt]
  5 & 2 & x & \shortstack[l]{\texttt{-0.669432458627258088807756766}\\\texttt{560685908183687401177881721303}\\\texttt{933846270762728020747989887159}\\\texttt{266}} \\[5pt]
  6 & 4 & y & \shortstack[l]{\texttt{0.4312821800183564134118023328}\\\texttt{329011468365299775114208035400}\\\texttt{405132844231227893005566564836}\\\texttt{18}} \\[5pt]
  7 & 3 & y & \shortstack[l]{\texttt{-0.593418257565607794147131052}\\\texttt{902379303831456603164876638088}\\\texttt{887598599913746499712706176082}\\\texttt{442}} \\[5pt]
  8 & 0 & y & \shortstack[l]{\texttt{-0.479012538959384353133079431}\\\texttt{525882160942433390265264451063}\\\texttt{351341069309509144360453723199}\\\texttt{162}} \\[5pt]
  9 & 7 & z & \shortstack[l]{\texttt{0.3950821090504650376521986546}\\\texttt{817323348681266067837696855183}\\\texttt{108328988934717397825504393081}\\\texttt{54}} \\[5pt]
  10 & 0 & z & \shortstack[l]{\texttt{0.4165153022927157118457819654}\\\texttt{650720972590496347096941588159}\\\texttt{771600622636956051128270647964}\\\texttt{37}} \\[5pt]
  11 & 1 & x & \shortstack[l]{\texttt{-0.296412601317025699166251982}\\\texttt{506222872498374630820152908105}\\\texttt{440866936332744700340974080221}\\\texttt{649}} \\[5pt]
  \bottomrule
\end{longtable}
\normalsize

\subsection{Middle Fixed Parameters}
\begin{center}\scriptsize
  \begin{tabular}{rrrr}\toprule
    $i$ & $x_i$ & $y_i$ & $z_i$\\\midrule
    0 & $\ast$ & -0.12508782251734035 & $\ast$ \\
    1 & $\ast$ & 0.19922464011814392 & $\ast$ \\
    2 & $\ast$ & $\ast$ & -0.3287408461854804 \\
    3 & $\ast$ & 0.3271508457133067 & -0.2448229081159972 \\
  \bottomrule\end{tabular}
\end{center}
Coordinates for labels $4,5,6,7$ are always imposed by $X_{7-i}=JX_i$.

\Needspace{12\baselineskip}
\subsection{Middle Variable Centers}
\small
\begin{longtable}{rrcl}
  \toprule Index & Vertex & Coordinate & Exact rational center \\\midrule
  \endfirsthead
  \toprule Index & Vertex & Coordinate & Exact rational center \\\midrule
  \endhead
  1 & 0 & z & \shortstack[l]{\texttt{0.6823277567493534973996458762}\\\texttt{857982016092983370833470159820}\\\texttt{198951410130769647437047277251}\\\texttt{86}} \\[5pt]
  2 & 1 & x & \shortstack[l]{\texttt{-0.846223033661651578013021651}\\\texttt{475536394479665383911212871022}\\\texttt{102865520467780161960855508602}\\\texttt{959}} \\[5pt]
  3 & 2 & y & \shortstack[l]{\texttt{0.5708162887883814053051200022}\\\texttt{926526134714886384884265717655}\\\texttt{389948381025082376691552442266}\\\texttt{38}} \\[5pt]
  4 & 3 & x & \shortstack[l]{\texttt{0.6051259216623858467931631451}\\\texttt{743102314553896350648469465032}\\\texttt{572658012747048323911090858412}\\\texttt{07}} \\[5pt]
  5 & 1 & z & \shortstack[l]{\texttt{-0.171482650283874236149634941}\\\texttt{145404438436154675113611844862}\\\texttt{387772780773411914616443994768}\\\texttt{474}} \\[5pt]
  6 & 2 & x & \shortstack[l]{\texttt{-0.363177831715389579303517962}\\\texttt{676784087879899523631996236534}\\\texttt{297641576991194708184622978601}\\\texttt{503}} \\[5pt]
  7 & 0 & x & \shortstack[l]{\texttt{0.5469701121653101349990476897}\\\texttt{423834889828524460554686472121}\\\texttt{252453370873921504184910073081}\\\texttt{79}} \\[5pt]
  \bottomrule
\end{longtable}
\normalsize

The reader can verify these exact coordinate tables by running
\path{code/proofs/check_manuscript.py}, which regenerates them using
\path{code/proofs/expand_paper.py}.

\section{The Exact Port Invariants}\label{app:invariants}

Each of the following constants is the ratio of the displayed integers.

\paragraph{Invariant $\kappa_0$}
\begin{flushleft}\small\ttfamily numerator: 34135779003287488671618031019671782250800097936990548192009204978\\[2pt] denominator: 11658805249772055355781839024560878088399285547058005194893959877\end{flushleft}
\paragraph{Invariant $\kappa_1$}
\begin{flushleft}\small\ttfamily numerator: 197963993380848124644671352310051480156731403415585401119461280032\\[2pt] denominator: 23586716069250660070577843160813818279890200073220144873845740953\end{flushleft}
\paragraph{Invariant $\kappa_2$}
\begin{flushleft}\small\ttfamily numerator: 36424441485755323185595901433419242799997826322327022519575870169\\[2pt] denominator: 3814110189865220980762038965553634106477269377045615636003558781\end{flushleft}

The reader can verify these rational constants against the saved matching
invariants by running \path{code/proofs/check_manuscript.py}.

\section{Certified Rational Angle Enclosures}\label{app:angles}

Here, we expand on and prove the correctness of the routine \path{code/proofs/rigorous_angles.py}. Let $p$ be the working binary precision, $B = p + 64$, $S = 2^B$, and $N = 2\lceil B/4\rceil$. 

\subsection{Reduced Arctangents and \texorpdfstring{$\pi$}{pi}}

For a rational $t \in [0, 1/2]$, Taylor expansion gives  
\begin{equation}\label{eq:atanseries}
  \arctan t = \sum_{k=0}^{N-1}\frac{(-1)^kt^{2k+1}}{2k+1} + \rho_N(t), \qquad 0\leq\rho_N(t)\leq\frac{t^{2N+1}}{2N+1}.
\end{equation}
We can be more precise. Let $L_0=\lfloor St\rfloor$, $U_0=\lceil St\rceil$ and
\[  A=\lfloor L_0^2/S\rfloor,\qquad 
  D=\lceil U_0^2/S\rceil,\qquad
  L_{k+1}=\lfloor L_k A/S\rfloor,\qquad
  U_{k+1}=\lceil U_k D/S\rceil.\]
These are all nonnegative integers. Induction gives $L_k/S\leq t^{2k+1}\leq U_k/S$. Define $a_k=\lfloor L_k/(2k+1)\rfloor$ and $b_k=\lceil U_k/(2k+1)\rceil$. Then an enclosing interval is
\begin{align*}
  \underline a(t) &= S^{-1}\left(\sum_{\substack{0\leq k<N\\ k\ {\rm even}}}a_k -\sum_{\substack{0\leq k<N\\ k\ {\rm odd}}}b_k\right),\\
  \overline a(t)&=S^{-1}\left(\sum_{\substack{0\leq k<N\\ k\ {\rm even}}}b_k -\sum_{\substack{0\leq k<N\\ k\ {\rm odd}}}a_k + \left\lceil\frac{U_N}{2N+1}\right\rceil\right).
\end{align*}
The last term bounds the positive remainder in \eqref{eq:atanseries}. For $t=0$ the code returns $[0,0]$ directly, and for negative reduced arguments it returns $[-\overline a(-t),-\underline a(-t)]$. The constant $\pi$ is enclosed using Machin's identity: 
\begin{equation}\label{eq:machin}
  \pi=16\arctan(1/5)-4\arctan(1/239).
\end{equation}
For completeness, write $\alpha = \arctan(1/5)$. The double-angle formula for $\tan$ gives $\tan(2\alpha) = 5/12$ and $\tan(4\alpha) = 120/119$, and the difference formula gives $\tan(4\alpha - \pi/4) = 1/239$, whence Machin's identity follows. Using \eqref{eq:machin}, if the enclosing intervals for the arctangents are $[a, b]$ and $[c, d]$, the enclosing interval for $\pi$ is then $[16a - 4d, 16b - 4c]$.  

For any positive rational $t$, we can evaluate $\arctan t$ as
\begin{equation}\label{eq:atanreduce}
  \arctan t=\begin{cases}
    \arctan t,&0<t\leq1/2,\\
    \pi/4+\arctan((t-1)/(t+1)),&1/2<t\leq2,\\
    \pi/2-\arctan(1/t),&t>2.
  \end{cases}
\end{equation}

\subsection{The Upper Half-plane Rectangle}

\begin{lemma}\label{lem:angleinclusion}
  For closed intervals $Y = [y_-, y_+]$ and $X = [x_-, x_+]$ with $y_- > 0$, the routine \code{atan2_upper(Y, X)} returns an interval containing $\operatorname{atan2}(y, x)$  for all $y, x \in Y \times X$.
\end{lemma}
\begin{proof}
  First we consider rational point arguments with $y > 0$. If $x > 0$, evaluate \eqref{eq:atanreduce} at $t = y/x$, if $x < 0$, evaluate it at $t = y/|x|$ and subtract the resulting interval from the rational enclosure of $\pi$, and if $x = 0$, halve the rational $\pi$ interval. These are precisely the three branches of $\operatorname{atan2}$ with range $(0, \pi)$, so Appendix \ref{app:angles}'s rational operations give an enclosing interval in each case. 
  
  Now we local the extrema. On the upper half plane, we have 
  \[\frac{\partial \theta}{\partial x} = -\frac{y}{x^2 + y^2} < 0, \qquad \frac{\partial \theta}{\partial y} = \frac{x}{x^2 + y^2}.\]
  For a fixed $y$, the minimum is at $x_+$ and the maximum at $x_-$. On the $x_+$ edge, we see that the minimum is obtained at $y_+$ when $x_+ < 0$ and at $y_-$ when $x_+ \geq 0$. On the $x_-$ edge, the maximum is obtained at $y_-$ when $x_- < 0$ and at $y_+$ when $x_- \geq 0$. Thus, we can enclose the minimum and maximum values and take the lower and upper bounds, respectively, to enclose the full range.

  Finally, we store rational numbers exactly, avoiding floating-points, and we use outward rounding for arithmetic. The routine proved here is implemented in \path{code/proofs/rigorous_angles.py}, and the skeptical reader can independently verify the angle enclosures used in the finite certificate calculation by running \path{code/proofs/audit_angle_enclosures.py}.
\end{proof}

\section{Rational Bounds for the Geometric Conditions}\label{sec:dyadic}

For each block, the root box $Q_B$ allows each free parameter to vary by at most $\ClaimRootRadius$ from its rational center. It is convenient to use dyadic boxes, so let $D = 2^{20}$. For each of the 24 coordinate entires, take its rational lower and upper endpoints in the root box, and round the lower endpoint down to the nearest multiple of $D^{-1}$ and the upper endpoint up to the nearest multiple of $D^{-1}$. The product of these intervals contains $Q_B$. Exact interval arithmetic gives the following lower bounds:
\begin{center}\small
  \begin{tabular}{lcc}\toprule
    Quantity & Cap & Middle \\\midrule
    $1-\|X_i\|^2$ & $1/3$ & $1/5$\\
    Squared face cross-product norm & $1/20$ & $1/8$\\
    $B_V$ and $B_W$ & $1/10$ & $3/20$\\
    $B_VB_W-B_{VW}^2$ & $1/50$ & $1/50$\\
    Squared port cross-product norm & $9/10$ & $3/10$\\
    Signed shared-edge determinant & $1/1000$ & $1/100$\\
    Signed separating functional value & $99/100$ & $99/100$\\
    Nonboundary height $z$ & $1/100$ & ---\\
    Port normal Lorentz norm squared & --- & $1/5$\\
    Strict inward port-plane margin & --- & $3/10$\\
    Negative outward-normal Lorentz product & --- & $3/8$\\
    Ultraparallel Gram gap & --- & $1/10$\\\bottomrule
  \end{tabular}
\end{center}
As in \eqref{eq:ultra}, let $N_i = (d_i, n_i), q_i = \langle N_i, N_i\rangle_\LL,$ and $L = \langle N_0, N_1\rangle_\LL$. The last four rows are lower bounds for $q_i$, $d_i - \langle n_i, X_j\rangle$ for every vertex $j$ outside port $i$, $-L$, and $L^2 - q_0q_1$, respectively. The skeptical reader can run \path{code/exact_geometry/check_certificates.py} to reconstruct the claims from the meshes and verify them. Run \path{code/verify_paper.py} to check the displayed bounds and the denominator against the computed values.

\section{Verification and Figures}\label{app:files}

The files \path{code/proofs/gluing_matched_cap_certified.json} and \path{code/proofs/gluing_matched_middle_certified.json} contain the exact rational centers, port invariants, and separating functionals used in the proof. The candidates were found numerically from Schwartz's seed and refined before these rational data were fixed. Verification starts from the supplied data and requires no search.

Commands in this paper are run from the package root with \code{python3 -E -B}. After installing the requirements, the following command repeats every calculation, regenerates graphics, runs tests, and compiles the paper:
\begin{center}
  \code{python3 -E -B reproduce.py --output ../origami-results}
\end{center}
Ensure the output directory is new and outside the package before running. For only mathematical verification, add \code{--only proofs}. For individual claims, each proof cites the file used for verification. 

To draw another example, say genus 4, run 
\begin{center}
  \code{python3 -E -B code/hyperbolic_origami.py --genus 4 --output ../genus4.png}
\end{center}
See the README for camera settings and export options. 

\end{document}